\documentclass[10pt]{amsart}

\newcommand{\be}{\begin{enumerate}}
\newcommand{\ee}{\end{enumerate}}
\newcommand{\bd}{\begin{description}}
\newcommand{\ed}{\end{description}}
\newcommand{\bi}{\begin{itemize}}
\newcommand{\ei}{\end{itemize}}

\newtheorem{thm}{Theorem}[section] %

\newtheorem{example}[thm]{Example} %

\newcommand{\eqn}{\begin{eqnarray}}

\newcommand{\eeqn}{\end{eqnarray}}

\usepackage[dvips]{graphics,color}

\begin{document}

\title{$SU(n)$ and $U(n)$ Representations of Three-Manifolds with Boundary }

\date{}

\maketitle

\centerline{Sylvain E. Cappell \ and \ Edward Y. Miller}

\vspace{.2in}
\large{

\centerline{\textbf{Abstract}}

{\small

Let $M^3$ be a compact, connected, oriented  three-dimensional manifold with non-empty  boundary,
 $\partial M$.
 This paper obtains results on extending on extending flat vector bundles or equivalently
 representations of the fundamental
 group  from $S$, part of the boundary, to the whole manifold $M^3$.
 The proof uses the introduction, investigation, and complete computation
 up to sign of new numerical invariants $\lambda_{U(n)}(M^3,S)$ (respectively,
 $\lambda_{SU(n)}(M^3, S)$), where $S$ is a proper connected subsurface of
 $\partial M$ which is connected.
 These numerical invariants
   `count with multiplicities and signs'
the number of representations up to conjugacy of the fundamental group of $M$, $\pi_1(M)$,
to the unitary group  $U(n)$ (resp., the special unitary group $SU(n)$)
which when restricted to  $S$
are conjugate to a specified irreducible representation
of $\pi_1(S)$.
These invariants are inspired by the
work of Casson  on $SU(2)$ representations of closed manifolds \cite{Akbulut1}.
 All the invariants treated here are seen to be
independent of the choice of irreducible representation, $\phi$,  on the surface $S$ in the boundary of $M^3$, $\partial M^3$.

If the difference of the Euler characteristics
$T =\chi(S) - \chi(M)$ is non-negative,  a $(dim \ U(n)) \times  T$-dimensional cycle is produced
that carries information about the space of such $U(n)$ representations.
For $T=0$,  the above   integer invariant  results.
 For $T > 0$, under the assumption that $\phi$ sends each
boundary component of $S$ to the identity,
a list  of  invariants results which are expressed as
a homogeneous polynomial in many variables, $\Lambda_{U(n)}(M,S,o)$.
Here $o$ records the relevant choices of orientations needed.
 (Such  a pattern of polynomial invariants is  reminiscent of the   work
of Donaldson on 4-manifolds \cite{Donaldson1}.) The $SU(n)$ case is  treated in a parallel fashion,
giving other polynomial invariants for $T >0$.

For $T=0$, the resulting integer invariants for $U(n)$  (resp., $SU(n)$)  are explicitly computed up to sign in all cases.
For $T>0$, an example is given to show that the $U(n)$ (resp., $SU(n)$) invariants are non-zero
in some cases.

Some   applications to problems  of extending $U(n)$  (resp., $SU(n)$)
representations to the fundamental groups of   three-manifolds, in particular to
the case of  rational homology cobordisms, are   given in \S 1.
After a historical review in \S 2,   \S 3  explains and states the main results  on the invariants, \S 4
defines the invariants,
\S 5 proves the theorems, while \S 6 treats stabilization issues.
\S 7 summarizes that, for $T=0$, the case of the  numerical invariant, there is
a $U(n)$ (resp., $SU(n)$) gauge theoretic reformulation, analogous to Taube's $SU(2)$ gauge theoretic
reformulation \cite{Taubes1} of Casson's invariant but here using a Fredholm operator $Z/2$ spectral flow.

}

\section{Some Applications to extending representations}

Flat vector bundles over a manifold correspond to linear representations of the fundamental group, with
the correspondence given by holonomy.
Extending representations from the fundamental groups of (part of) the boundary
of a manifold to the fundamental group of the whole manifold, or equivalently, extending flat bundles, are important in defining
and investigating invariants in many contexts in both high and low dimensional topology and geometry.
This section states results on this problem of extending representations of the fundamental groups
from (part of) the boundary of a three-manifold. These results are  obtained using  the invariants defined and
investigated in later sections. Theorem \ref{theoremfirst} treats the case of connected
boundary and theorem  \ref{theoremapplication} the case when the boundary consists
of two boundary components and the three-manifold is a rational homology cobordism
between them.

\begin{thm}  \label{theoremfirst} Let $M^3$ be an oriented, connected, compact three-manifold with connected boundary
and $S$  a connected subsurface of its boundary  for which
the inclusion  of $S$ in $M^3$, $i : S \subset M^3$, induces an isomorphism of first rational homology groups, i.e.,
$$
i_\star : H_1(S,Q) \stackrel{\cong}{\rightarrow} H_1(M^3,Q),
$$
Then  picking a base point $p$ in $S$:

\begin{itemize}
\item[1.] Case $U(n)$: For any $U(n)$ representation $\phi: \pi_1(S,p) \rightarrow U(n)$
there exists a representation $\phi' : \pi_1(M^3,p) \rightarrow U(n)$ which restricts on $S$
to $\phi$.
\item[2.] Case $SU(n)$: If moreover, $\phi$ is irreducible and has image in $SU(n)$, then its
extension $\phi'$ may be chosen to enjoy the same  property.
\end{itemize}
\end{thm}

This result covers many topologically quite different possible choices
 of the  surface $S$ which may have one or many boundary circles. For example,
if $M^3$ is a three-manifold for which the inclusion of the boudary, a Riemann surface of genus $g$,
$\partial M \subset M^3$ induces a surjection $H_1(\partial M,Q) \rightarrow H_1(M^3,Q)$
of rational homology,  $S$ may be chosen to  be a thickening in $\partial M$ of a wedge of
circles suitably chosen, or if $g$ is even, a suitable once-punctured Riemann surface
of genus $g/2$, or many other possibilities satisfying the homological condition  of theorem
\ref{theoremfirst}.

\vspace{.2in}

A three-dimensional, oriented, \textit{rational homology cobordism} consists of
a compact, connected, oriented 3-manifold $W$ with two boundary components, denoted
$N_1,N_2$, such that the maps induced by inclusion on first and second rational homology
$$
H_1(N_1,Q) \rightarrow H_1(W,Q) \ and \ H_2(N_1,Q) \rightarrow H_2(W,Q)
$$
are isomorphisms. By Poincar\'{e} duality, this is
equivalent to the same condition on the other boundary component,
$$
H_1(N_2,Q) \rightarrow H_1(W,Q) \ and \ H_2(N_2,Q) \rightarrow H_2(W,Q)
$$
are isomorphisms; or again equivalently, to the condition that both
$H_1(N_1,Q) \rightarrow H_1(W,Q)$ and $H_1(N_2, Q) \rightarrow H_1(W,Q)$
are  isomorphisms.

It is of natural topological  interest to understand to what extent  a homology cobordism
behaves like a cylinder.
 A theorem of this nature  is proved in this paper :

\begin{thm} \label{theoremapplication}

Let $W^3$ be a compact, connected, oriented
3-manifold with  boundary
consisting of two connected components, denoted $N_1,N_2$, such that
the inclusion $i : N_1 \subset W^3$ induces an isomorphism
of rational  first homology,
$$
i_\star : H_1(N_1,Q) \stackrel{\cong}{\rightarrow}  H_1(W^3,Q).
$$
Then  $W^3$ is a rational homology cobordism and, picking a base point $p$ in $S$:

\begin{enumerate}

\item[1.] Case $U(n)$: For any $U(n)$ representation,
$$
\phi : \pi_1(N_1,p) \rightarrow U(n),
$$
there exists a representation $\phi':\pi_1(W^3,p) \rightarrow U(n)$
which restricts on $N_1$ to  $\phi$.

\item[2.] Case $SU(n)$: If moreover, $\phi$ is irreducible and has image in $SU(n)$,
then its extension $\phi'$ may be chosen to enjoy the same property.
\end{enumerate}

\end{thm}

\section{Historical background on representation counting} \label{sectionHistorical}

Let  $M^3$ be  a closed,  connected, oriented 3-manifold which is
a homology sphere, that is, its integral homology is the same as that of the
standard 3-sphere,  $H_\star(M^3,Z) \cong H_\star(S^3,Z)$.  Casson
introduced an  invariant which `counts with signs and multiplicities'
the number of representations up to conjugacy of the fundamental group of $M^3$,\ $ \pi_1(M^3)$,
to the special unitary group $SU(2)$ up to conjugacy
[possibly after suitable  perturbation] \cite{Akbulut1,Saveliev1, Saveliev2}.
 He further derived surgery formulae for effectively computing his invariant
and established relations to the Rohlin invariant,  and  other subjects.
Allied invariants have since been studied by Boyer and Nicas \cite{BoyerNicas1}, Walker \cite{Walker1}
[rational homology spheres], Lescop \cite{Lescop1} [general closed oriented three-manifolds],  Lin \cite{Lin1} [knot invariants],
Herald \cite{Herald1} [knot invaiants], Frohman and Long \cite{FrohmanLong1}  [ knot invariants],
 Ruberman and Saveliev \cite{RubermanSaveliev1} [four dimensional invariants].

In a different direction, Casson's
invariant was reinterpreted in an analytical gauge theoretic fashion by
Taubes \cite{Taubes1}, forging an important further
  bond between
geometric  topology and  mathematical physics.
The topological sign assigned by Casson to a representation that
is transverse is in Taubes' setting  expressed in terms of a spectral
flow. The works of Casson and Taubes have been the inspiration for much subsequent
work;  a notable example was the work of Floer \cite{Floer1,Floer2, Floer3} who
gave an instanton theory refinement and a $Z/8$  graded homology theory with the grading
given by a  spectral flow modulo 8. The Euler characteristic of this Floer homology
is the numerical $SU(2)$ invariant of Casson and Taubes  for three-dimensional homology
spheres.

Because of  difficulties in treating the singular strata which correspond to reducible representations
in the spaces of representations (see \cite{Goldman1,Goldman2})
to $SU(n)$ or $U(n)$, Casson's original
invariant for closed 3-manifolds was restricted to the group $SU(2)$ and to integral homology 3-spheres.
In the decades since, substantial efforts expanded the class of 3-manifolds which can be treated, but
have only very  partially lifted the restriction on the Lie group. The difficulties
escalate as $n$ increases, as the number and kinds of singular strata increase in the associated
symplectic varieties of $SU(n)$ ( resp., $U(n)$) representations.

For $G = SU(2)$, a definitive picture emerged from a series of extensions. The work of
Boyer and Lines \cite{BoyerLines1}  extended the invariant to the case where
 the fundamental group  equals $Z/n$ by adding rational corrections for each reducible representation.
 The work of
Walker \cite{Walker1} completely extended the  invariant  to
rational homology spheres, $H_\star(M^3,Q) \cong H_\star(S^3,Q)$, defining an invariant
valued in the rationals $Q$ and  proving
all the analogues of Casson's original treatment. The calculations of Boyer and Lines
 involve Dedekind sums associated to the reducible representations of lens spaces,
and in Walker's work their generalizations. The work of Lescop \cite{Lescop1}
employing a topological-combinatorial approach treated in generality  oriented closed 3-manifolds,
still for $G = SU(2)$. The work of Boden and Nicas  treats an $SU(n)$ knot invariant for $n \ge 2$
\cite{BodenNicas1}, as does the work of Frohman \cite{Frohman1}.

Various extensions to the case of $SU(3)$ have been pursued by several authors
in the  standard topological  Casson setting and in the gauge setting
of Taubes.  For example, there is an  $SU(3)$ invariant of  Boden and Herald \cite{BodenHerald1}
for integral homology 3-spheres with  explicit  computations  carried out
  by  Boden, Herald, Kirk, Klassen \cite{ BodenHeraldKirkKlassen1,BodenHeraldKirk1}.
A different $SU(3)$ Casson type invariant was introduced by Ronnie Lee
and the present authors; this invariant is perturbative  \cite{CappellLeeMiller1}.
However, to date the increasing difficulties with the stratification
of the representation spaces have precluded a general definition
 for  $SU(n)$ representations, for $n > 3$,
for closed 3-manifolds.

\vspace{.2in}
As seen in this paper, by considering representations to $U(n)$ (respectively, $SU(n)$) which restrict
to a specified irreducible representation  to $U(n)$ (resp., $SU(n)$) on part of the boundary of a
non-closed manifold, the, in general, as yet unresolved,  technical difficulties
caused by the need for the closed case  to treat the contributions of the
singular strata are absent. As seen below,  the $U(n)$ (resp., $SU(n)$) numerical  analogue of the Casson invariant
in the present context
can then be computed homologically for all $n$, at least up to sign.

\section{Results on numerical and polynomial invariants} \label{sectionStatements}

The present section  deals with Casson type invariants in  the context
of  compact, connected, oriented  3-manifolds $M^3$ with non-empty  but connected boundary $\partial M^3$
together with a specified embedding
$$
S_1 \subset \partial M^3
$$
of a connected [non-closed]  subsurface $S_1$ into   $\partial M^3$ with base point $p$,
and a choice of \textbf{irreducible} representation of the fundamental group of $S_1$,
$$
\phi : \pi_1(S_1,p) \rightarrow U(n),
$$
to the  unitary group $U(n)$. The underlying idea is to
construct invariants that measure the number of representations
of the fundamental group of $M^3$, $ \pi_1(M^3)$, to $U(n)$
which restrict on $S_1$ to  the chosen representation, $\phi$,
up to conjugation. It is then shown
 that the resulting invariants are independent of the choice of such  irreducible $\phi$.

\vspace{.2in}
An entirely  parallel analysis for the special unitary group $SU(n)$
proceeds in a like manner with minor changes and is indicated below.
More generally, there are also definitions and results analogous to the results of this section for
counting irreducible representations to any compact, connected Lie group $G$. For such general
Lie groups irreducible means, as usual,  that the only elements of $G$ which  commute with all
elements of the image of the given representation are those in the center of the group $G$.

\vspace{.2in}
The initial restriction to \textit{irreducible} representations $\phi$ is  crucial
in this section; however,
since any representation to $U(n)$ is a sum of irreducible representations,
general results may often be treated by first studying this special case, as in the first statement of
theorems \ref{theoremfirst} and \ref{theoremapplication}.
 Moreover,  the case of manifolds
with several boundary components may often be reduced to the case of one, by connecting the boundary components
by disjointly embedded paths and deleting a tubular neighborhood of them.
By these means the results on this special case have general implications,
for example theorem \ref{theoremapplication}.

\vspace{.2in}
The aim is  to ``count appropriately  with signs and multiplicities''
the number of
representations, up to conjugacy,  of the fundamental group of $M^3$, $\pi_1(M^3)$, to $U(n)$
(respectively, $SU(n)$)
which restrict to $\phi$, up to conjugacy. To deal with issues of global sign, additionally a
choice of orientations of $M^3$ and of the real cohomology groups   $ H^1(M,R), H^2(M,R), H^1(S_1,R)$ is needed, see
\S \ref{sectionstable}. These
orientation choices
are recorded via the symbol $o$.

\vspace{.2in}
Firstly, it will be shown that this can
be carried out in the case that there is equality of Euler characteristics,
$$
\chi(S_1) = \chi(M^3),
$$
to yield  integer valued invariants,
$$
\begin{array}{c}
\lambda_{U(n)}(M^3,S_1) \in Z, \\
\lambda_{SU(n)}(M^3,S_1) \in Z, \\
\end{array}
$$
which are  shown to be independent of the irreducible representation $\phi$ chosen on $S_1$.
The
absolute values $|\lambda_{U(n)}(M^3,S_1)|, |\lambda_{SU(n)}(M^3,S_1)|$  are  independent
of the orientation choices  recorded by $o$. Let $p$ be the base point of $S_1$
chosen to lie in the boundary of $S_1$.

\vspace{.2in}
In contradistinction to the original $SU(2)$ Casson invariant, these invariants
$\lambda_{U(n)}(M^3,S_1) \in Z $ and $  \lambda_{SU(n)}(M^3,S_1) \in Z$
will be shown to be of a  homological character, as seen from the following theorem
which, in particular, completely determines them up to sign.

\begin{thm}{Case $U(n)$:} \label{theoreminv1} In the case that $\chi(S_1) = \chi(M^3)$ and $\partial M^3$ connected:

\begin{enumerate}
\item If $\lambda_{U(n)}(M^3,S_1) \neq 0$, then there is at least
one representation $\phi': \pi_1(M^3,p) \rightarrow U(n)$ restricting
to $\phi : \pi_1(S_1,p) \rightarrow U(n)$.
\item  The value of  the
invariant $\lambda_{U(n)}(M^3,S_1)$ is independent of the choice of
irreducible representation $\phi$.
\item If the rational homology $H_2(M^3,Q)$ is non-vanishing or the induced mapping on rational homology
$$
H_1(S_1,Q) \rightarrow H_1(M^3,Q)
$$
is not an isomorphism, then $\lambda_{U(n)}(M^3,S_1) = 0$.
\item If  the induced mapping on rational cohomology
$$
H_1(S_1,Q) \stackrel{\cong}{\rightarrow} H_1(M^3,Q)
$$
is  an isomorphism, then  $H_2(M_3,Q)= 0$ and  the absolute value $|\lambda_{U(n)}(M^3,S_1)|$
equals \newline
$K^{n}$, where $K$ is the order of the finite
abelian group \newline $H^2(M^3,S_1,Z)$, the second cohomology. In particular,
$\lambda_{U(n)}(M^3,S_1) \neq 0$ in this case.

\item $ \lambda_{U(n)}(M^3,S_1)=(-1)^n \ \lambda_{U(n)}(-M^3, S_1)$ where $-M^3$ is
$M^3$ with the opposite orientation.

\end{enumerate}
\end{thm}

\begin{thm}{Case $SU(n)$:}  \label{theoreminv1SU}
The same results hold replacing $U(n)$ by $SU(n)$ except that in statement (4), $K^n$ is replaced
by $K^{n-1}$ and in statement (5), $(-1)^n$ is replaced by $(-1)^{n-1}$.
\end{thm}

These changes from $K^n$ and $(-1)^n$ to $K^{n-1}$ and $  (-1)^{n-1}$ reflect that the rank of $U(n)$
as a Lie group  is
$n$ and the rank of  $SU(n)$ is $n-1$.

\vspace{.2in}

Recall that for an arbitrary compact, oriented, odd dimensional manifold $N$
with boundary $\partial N$, their Euler characteristics are related by:
$$
\chi(N) = (1/2) \ \chi(\partial N).
$$
This is readily verified  by forming the double of $N$, called $\tilde{N}$, obtained
by taking two disjoint copies of $N$ with the tautological identification
on their boundaries. By Poincar\'{e} duality for the closed odd dimensional manifold $\tilde{N}$, the Euler
 characteristic of $\tilde{N}$ vanishes and so, by additivity of the Euler characteristic,
one obtains  $0 = \chi(\tilde{N}) = 2 \ \chi(N) - \chi(\partial N)$.

In the present context, if  the connected surface $\partial M^3$ has genus $g$,  the Euler characteristic of $M^3$
is just
$$
\chi(M^3) = (1/2) \ \chi(\partial M^3) = (1/2)(2-2g) = 1-g.
$$
So the condition that $\chi(M^3) = \chi(S_1)$ of theorem \ref{theoreminv1}  is just that
$\chi(S_1) = 1-g$. That is, $H_1(S_1,Z) \cong Z^g$.
For example, $S_1$ could be under the Euler characteristic condition be chosen to
be a regular neighborhood of a bouquet of $g$ circles imbedded into $\partial M^3$
which has genus $g$; but as noted in the discussion after theorem \ref{theoremfirst}
there are in general many other topologically distinct possible choices for $S_1$.

\vspace{.2in}

Secondly and correspondingly, if $\chi(S_1) - \chi(M^3) = T  >0$ is even, [note, $T = (g - g_1)$ with $\chi(S_1) =1-  g_1$],
 then (in analogue to Donaldson's
invariants of four dimensional manifolds \cite{Donaldson1})  additional invariants will be defined,
under the mild assumption that $\phi$ sends each homotopy class representing a boundary
component of $S_1$ to the identity in $U(n)$.

Note : This additional assumption excludes a few low genus cases. More precisely,
let  the connected surface $S_1$ be obtained from a closed Riemann
surface of genus $k$  by deleting $l$ 2-disks. In particular,
$\chi(S_1) = 1-g_1 = (2 - 2k)  -l = 1 -( 2k +l-1)$ or $g_1 = 2k+l-1$ and $S_1$ has
$l$ boundary components.  A representation
with the above property will factor through the fundamental
group of the closed surface. In particular, it will be
reducible if $k=0,1$ for $n \ge 2$.
 Thus the discussion of polynomial invariants below implicitly uses the condition
that $g_1 $ is greater than the number of components of the boundary of $S_1$ plus one.
With this assumption, any irreducible representation of $\pi_1(S_1)$ can be deformed
through irreducible representations to an irreducible representation
that sends each boundary component to the identity.

Explicitly, let $I,J$ (possibly vacuous) be two   multi-indices of pairs of non-negative integers,
$I = [(i_1,r_1), (i_2,r_2), \cdots, (i_a,r_a)], \ J =  [(j_1,s_1), (j_2,s_2), \cdots, (j_b,s_b)] $ such that
$i_1 < i_2< \cdots i_a$ and $ j_1 < j_2 < \cdots < j_b$ with
$$
 T = \chi(S_1)- \chi(M^3) = ( \ \Sigma_{p=1}^a ( 2 \ i_p)(r_p)  \ + \ \Sigma_{p=1}^b (4 \ j_p-2)(s_p).
$$
Then one   defines a  rational invariant
$$
\lambda_{I,J,U(n)}(M^3, S_1,o)
$$
 which
``counts with signs and multiplicities'' the
number of the $(dim \ U(n)) \times T$ parameter families of representations which extend to $M^3$
the irreducible representation $\phi$, where $\phi$ is chosen to send each boundary component
of $S_1$ to the identity element $Id \in  U(N)$.  As explained in \S \ref{sectionDefs1}, $o$
records a choice of the relevant orientations those of $M^3$, and  $H^2(M,R), H^1(M,R), H^1(S_1,R)$.

All these may be recorded by a single homogeneous polynomial invariant,
$$
\Lambda_{U(n)}(M,S_1,o) = \Sigma_{I,J} \ \lambda_{I,J,U(n)}(M,S_1,o) \ X_{I,J},
$$
with the sum over multi-indices $I,J$, as above and
$X_{I,J} $ denoting  the monomial :
$$
X _{I,J}= ( x_{i_1}^{r_1} x_{i_2}^{r_2} \cdots x_{i_a}^{r_a})  ( y_{j_1}^{s_1} y_{j_2}^{s_2} \cdots y_{j_b}^{s_b})
 $$
Here one should regard $x_i$ to be in degree $2  i$ and $y_j$ to be in degree $4j-2$.

\vspace{.2in}
The corresponding multi-index invariants for $SU(n)$, $\lambda_{I,J,SU(n)}(M^3,S_1,o)$,  are indexed by multi-indices
(possibly vacuous), $ I = [(i_1,r_1), (i_2,r_2), \cdots, (i_a,r_a)], \ J =  [(j_1,s_1), (j_2,s_2), \cdots, (j_b,s_b)] $, of non-negative
indices with  $i_1 < i_2< \cdots i_a$ and $ j_1 < j_2 < \cdots < j_b$ subject to the added constraints $i_1 >1$, if $I$ is present,
and
$ j_1> 1$, if $j_1$ is present. All these may be recorded by a single homogeneous polynomial invariant
$$
\Lambda_{SU(n)}(M,S_1,o) = \Sigma_{I,J} \ \lambda_{I,J,U(n)}(M,S_1,o) \ X_{I,J}
$$
with the sum over multi-indices $I,J$, as above with the additional  constraints $i_1 >1$, if $I$ present and  $j_1>1$,
if $J$ is present.

\begin{thm}{Case $U(n)$ : } \label{theoreminv2}
Let $\phi$ send each  homotopy  classes representing a  boundary component
of $S_1$ to the identity in $U(n)$. Then
with
 $T = \chi(S_1)- \chi(M^3)$ one has:
\begin{enumerate}
\item If $\Lambda_{U(n)}(M^3,S_1,o) \neq 0$, then there is a
$( dim \ U(n)) \times  T$ dimensional family
of  representations $\{\phi'_\lambda \}: \pi_1(M^3) \rightarrow U(n)$,  each restricting
to $\phi : \pi_1(S_1) \rightarrow U(n)$.
\item  The value of  the
invariant $\Lambda_{U(n)}(M^3,S_1,o)$ is independent of the choice of
irreducible representation $\phi$ sending  the homotopy
class of each boundary component of $S$  to the identity.

\item The invariant $\Lambda_{U(n)}(M^3,S_1,o)$  changes by the sign $(-1)^{n}$ under
a change of orientation of any of  $M^3$, $H^1(M,R),H^2(M,R), H^1(S_1,R)$.
In particular,  the absolute value of the coefficient of any monomial term in  the polynomial
$\Lambda_{U(n)}(M,S_1,o)$ is  independent
of the orientation choices  $o$.
\end{enumerate}
\end{thm}

\begin{thm}{Case $SU(n)$ : } \label{theoreminv2SU}
The same results, but with the added constraint on indices  $i_1>1$, if $I$ is present,
and $j_1 >1$, if $J$ is present,
   hold replacing $U(n)$ by $SU(n)$ except that in statement (4), $K^n$ is replaced
by $K^{n-1}$ and in statement (5),  $(-1)^n$ is replaced by $(-1)^{n-1}$.
\end{thm}

It would be interesting to have formulae for the polynomial invariants $\Lambda_{U(n)}(M,S_1,o),\newline
\Lambda_{SU(n)}(M,S_1,o), $
analogous to theorems \ref{theoreminv1}, \ref{theoreminv1SU},  parts 3 and 4.
An example is given at the end of section \ref{sectionstable} to show that
these invariants do not always vanish.

There are for general compact, connected, Lie groups $G$ definitions and
results analogous to the results of this section for `counting with signs
and multiplicities' the number of representations which restrict to a conjugate
of a prescribed, irreducible representation on part of the  boundary; and
more generally , there are polynomial invariants associated to $G$ when $T>0$.
These are all again independent of the  choice of the irreducible  representation
prescribed on part of the  boundary.

 The present treatment for 3-manifolds
with boundary of some representation-theoretic invariants   which had encountered symplectic and topological obstacles
for closed
three-manifolds, may suggest a paradigm in some other settings, e.g., for other invariants
of three-manifolds with boundary.

\section{Definition of the Invariants } \label{sectionDefs1}

Let $M^3$ be a compact, connected, oriented  3-manifold with  non-empty,
connected boundary, $\partial M^3$, of genus $g$. Let
$$
\partial M^3 = S_1 \cup S_2 \ with \ \partial S_1 = \partial S_2 = S_1 \cap S_2
$$
be a decomposition of the boundary of $M^3$ into two non-empty connected subsurfaces, $S_1,  S_2$.
Pick a point  $p \in \partial S_1 = S_1 \cap S_2$ to serve as the base point for all fundamental  groups.

It is assumed   that there is chosen a  fixed irreducible representation, called $\phi$,
$$
\phi : \pi_1(S_1) \rightarrow U(n)
$$
of the fundamental group of $S_1$ to $U(n)$. The treatment of $SU(n)$ is
entirely parallel, replacing $U(n)$ by $SU(n)$, except where noted in the
computations.

Following Casson's approach \cite{Akbulut1},  choose handlebodies $H_1,H_2$ of
genus $h_1,h_2$, respectively, such that $M^3$ is a union of $H_1,H_2$ adapted to the
decomposition of $\partial M^3$. That is, firstly $H_1$ is obtained from the
regular neighborhood $S_1 \times [0,1]$ in $M^3$ of $S_1 = S_1 \times 0$ by adding
one handles; secondly  $H_2$ is obtained from the
regular neighborhood $S_2 \times [0,1]$ in $M^3$  of $S_2 = S_2 \times 0$ by adding
one handles. Hence, the fundamental groups of $ \pi_1(H_1),\pi_1(H_2)$ are free groups
obtained from the fundamental groups $\pi_1(S_1), \pi_2(S_2)$, also free groups,
by adding added generators, $h_1-g_1, h_2-g_2$, of them,
respectively. These handlebodies, $H_1,H_2$, are chosen to intersect
in a subsurface $U \subset M^3$ with the added  properties:
$$
\begin{array}{l}
M^3 = H_1 \cup H_2 \ with \ S_1 = H_1 \cap (\partial M^3) \ and \ S_2 = H_2 \cap ( \partial M^3) \\
U = H_1 \cap H_2 \ with \ \partial H_1 = S_1 \cup U \ and \ \partial H_2 = S_2 \cup U \\
U \ meets \ \partial M^3 \   transversally \ along \ S_1 \cap S_2 \\
U \ is \ connected
\end{array}
$$
Recall that a handlebody of genus $h$ is a regular neighborhood of a bouquet
of $h$ circles embedded in $R^3$; equivalently,  it is the union of a closed 3-ball, $D^3$,  together
with $h$ 1-handles, $ D^2 \times [0,1]$ suitably attached along the boundaries
pieces $D^2 \times 0, D^2 \times 1$ to disjoint parts of the boundary $\partial D^3$. In particular,
if the handlebodies, $H_1,H_2$ have genus $h_1,h_2$, then $\chi(H_1) = 1-h_1,\chi(H_2) = 1- h_2$
and the fundamental groups, $\pi_1(H_1),\pi_1(H_2)$ are free groups on $h_1,h_2$ generators,
respectively.

\vspace{.2in}
The following notation is used: For a space $X$ and subspace $S \subset X $  with  base point $p \in S$, let $R^\#(X)$ denote the space of representations
of the fundamental group, $\pi_1(X,p)$, to $U(n)$ and $R(X)$ be the quotient of $R^\#(X)$ by the conjugation
action of $U(n)$. For a specified representation
$\phi : \pi_1(S,p) \rightarrow U(n)$, define
 $R^\#(X,S,\phi)$  to be the space of representations of $\pi_1(X,p)$ to $U(n)$
which when restricted to $\pi_1(S,p)$ give the specified representation $\phi$.
Let $R(X,S,[\phi])$ be the image in $R(X)$ of  $R^\#(X,S,\phi)$ under the natural projection.
Alternatively,  $R(X,S,[\phi])$ may be described as the space of representations
of $\pi_1(X,p)$ to $U(n)$  which on $S$ equals $\phi$ up to conjugation, divided out by the conjugation
quotient action. In particular, it only depends on $\phi$ up to conjugation.

\vspace{.2in}

Note  that since $\phi : \pi_1(S_1,p) \rightarrow U(n)$ is assumed to be \textbf{irreducible}, the
natural quotient mapping
$
R^\#(M^3,S_1,\phi) \rightarrow R(M^3,S_1,[\phi])
$
sending a representation $\rho \in R^\#(M^3)$ to its image in  $R(M^3)= R^\#(M^3)/conjugation$
\textbf{ is a bijection}:
$$
R^\#(M^3,S_1,\phi) \cong R(M^3,S_1,[\phi])
$$
This holds since the irreducibility of $\phi$ means
that the conjugation action of $U(n)$ is free modulo the action of the
center $Z(U(n)) \cong U(1)$ of $U(n)$ which conjugates trivially. Otherwise expressed,
the projective unitary group $PU(n) = U(n)/Z(U(n))$ acts freely on the irreducible representations.

Hence, by $\phi$ irreducible, to suitably count with signs and multiplicities
the number of points in $R(M^3,S_1,[\phi])$, it suffices to consider
the space of representations $R^\#(M^3,S_1,\phi)$ instead, a useful simplification.

\vspace{.2in}

To give an  explicit model of $H_1,H_2,U = H_1 \cap H_2$ take a triangulation, say $K$ of $M^3$ for which
 $S_1 = |L_1|, S_2 = |L_2|$, for subcomplexes $L_1,L_2$ of $K$.  Now in the second baricentric subdivision, $K''$, let $H_a$
be the union of  closed simplicies of $K''$ which contain as a vertices of $K''$
the baricenter of a vertex of $K$, i.e., that vertex again, or the baricenter of
an edge of $K$.  Correspondingly, let $H_b$ be the subcomplex of $K''$ which is
the union of the closed simplicies of $K''$ which contain as a vertex of $K''$  a baricenter  of a 3-simplex of $K$
or a baricenter of a 2-simplex of $K$.
It is a standard fact that $H_a, H_b$ are a handlebodies with $M^3 = H_a \cup H_b$
and $H_a \cap H_b$ a separating surface  for $M^3 $\cite{Singer1}. Such a decomposition is
called a Heegaard decomposition.  Now let $H'_1$ be obtained by adjoining to $H_a$
those closed simplicies of $K''$ which have a vertex which is a baricenter of
a 2-simplex of the subcomplex $S_1$. Correspondingly,  let $H'_2$ be obtained by adjoining to $H_a$
those closed simplicies of $K''$ which have a vertex which is a baricenter of
a 2-simplex of the subcomplex $S_2$. This gives a handlebody decomposition of $M^3$ with
the property that there is an isotopy of $M^3$ carrying  the intersections
$H'_1 \cap \partial M^3 , H'_2 \cap \partial M^3 $ to $S_1,S_2$ respectively. Let $H_1,H_2$
be the respective images of $H_1',H_2'$. This constructs, as desired,  handlebodies adapted
to the decompostion $\partial M^3 = S_1 \cup S_2$.

 In this model, the pieces, $H_1,H_2,S_1,U = H_1\cap H_2, $
are all connected, each is homotopy equivalent to a bouquet of circles, and their fundamental groups are free groups.
Also, the fundamental group of $H_1$ is the free group on $h_1$ generators obtained by adding $h_1-g_1$ free generators
to the fundamental group of $\pi_1(S_1)$; similarly, the fundamental group of $H_2$ is the free group on $h_2$ generators
obtained by adding $h_1-g_1$ free generators
to the fundamental group of $\pi_1(S_1)$.

\vspace{,2in}

In this model, $H_1$ is the union of a regular neighborhood, $S_1 \times [0,1]$, of $S_1$ in $H_1$
and added 1-handles, say $j$ of them.
By isotopying these added 1-handles,  they may be arranged to
 attach to a fixed 2-disk, say $D^2 \subset S_1 \times 1 \subset M^3$,
That is, $H_1$ is exhibited as the boundary connected sum
of $S_1 \times [0,1]$, a handlebody of genus $g_1$, and a handlebody
$\hat{H}$ of genus $j$, with common intersection the 2-disk, $D^2$. Thus, $h_1 = g_1 + j$.
In particular, there is are homeomorphisms
$$
\begin{array}{l}
R^\#(S_1) \cong U(n)^{g_1},  \ R^\#(H_1) \cong U(n)^{h_1}, \   and  \ R^\#(H_1,S_1,\phi) \cong U(n)^{h_1-g_1} \\
and \  in \ a \ parallel \ manner  \\
R^\#(S_2) = U(n)^{g_2},   \ R^\#(H_2) \cong U(n)^{h_2},
\end{array}
$$
since the inclusion $\pi_1(S_1) \subset \pi_1(H_1)$ merely adds $j = (h_1- g_1)$ new free
generators. A representation in $R^\#(H_1,\phi)$ will be specified on the
$g_1$ free generators of $\pi_1(S_1)$ by $\phi$,  but arbitrary on the remaining $j = h_1-g_1$
free generators. With these identifications $R^\#(H_1) = U(n)^{g_1} \times U(n)^{(h_1-g_1)} $,
also.

Expressed in another fashion, there are homeomorphisms,
 $$
 \begin{array}{l}
 R^\#(H_1) \cong U(n)^{h_1}, \ R^\#(H_2) \cong U(n)^{h_2}, \\
 R^\#(S_1) \cong U(n)^{g_1}, \  R^\#(S_2) \cong U(n)^{g_2}, \
  and \  R^\#(U) \cong U(n)^u \\
  \end{array}
 $$
 where $\chi(H_1) = 1- h_1, \chi(H_2) = 1- h_2$ and $\chi(S_1) = 1-g_1, \ \chi(U) = 1-u$.

Moreover, the inclusions $S_1 \subset H_1, S_2 \subset H_2$ induce the surjections
$R^\#(H_1) = U(n)^{h_1} \rightarrow R^\#(S_1) = U(n)^{g_1} $ and
$R^\#(H_2) = U(n)^{h_2} \rightarrow R^\#(S_2) = U(n)^{g_2} $
defined by projecting to the first $g_1,g_2$ factors, respectively.

\vspace{.2in}

Next, note that $M^3 = H_1 \cup H_2$  with $U = H_1 \cap H_2$,
so the fundamental group $\pi_1(M^3,p)$ is, by van Kampen's theorem,  the amalgamated free product
of $\pi_1(H_1,p)$ and $\pi_1(H_2,p)$ along $\pi_1(U,p)$. Hence,
a representation $f : \pi_1(M^3,p) \rightarrow U(n)$ is uniquely specified
by its restrictions to $\pi_1(H_1,p)$ and $\pi_1(H_2,p)$;  and given a pair
of such representations, they arise from a representation of $\pi_1(M^3,p)$
if and only if  their restrictions to $\pi_1(U,p)$ are equal. That is,
as Casson had observed in the analogous closed  manifold setting,
$
R^\#(M^3) $ is precisely the intersection of the images of $R^\#(H_1)$ and $R^\#(H_2)$
in $R^\#(U)$.
$$
R^\#(M^3) =  [ Image \ R^\#(H_1) ] \cap [ Image \ R^\#(H_2) ]  \ in \ R^\#(U)
$$
In particular,
$$
\begin{array}{l}
R^\#(M^3,S_1, \phi)
=  [ Image \ R^\#(H_1,S_1,\phi) ] \cap [ Image \ R^\#(H_2) ]
 \ in \ R^\#(U) \\
 while \ \ R^\#(M^3,S_1, \phi) \cong R(M^3,S_1, [\phi])
\end{array}
$$

Let $K:R^\#(H_1) \rightarrow R^\#(U), \  L:R^\#(H_2) \rightarrow R^\#(U)$
be the mappings defined by the homomorphisms of  fundamental groups $\pi_1(U) \rightarrow  \pi_1(H_1), \ \pi_1(U) \rightarrow \pi_1(H_2)$
induced by the inclusions $U \subset H_1, U \subset H_2$.  Let $J : R^\#(H_1,\phi ) \subset R^\#(H_1)$
be the inclusion. In this notation, by $\phi$ irreducible, the set
$$
 R(M^3,S_1, [\phi]) \cong R^\#(M^3,S_1, \phi)
$$
to be counted ``with signs and multiplcities'' is identified
with the image via $(a,b) \mapsto L (b)$ of the set of pairs \newline
$\{ (a,b) \ | \ a \in R^\#(H_1,\phi), \ b\in R^\#(H_2) \
with \ K \cdot J(a) = L(b) \}$.

\vspace{.2in}
Now, the fundamental group of the connected oriented  surface $U$ is a free group on, say,  $u$ generators, so picking
standard generators
 gives an identification
 $$
 R^\#(U) \cong U(n)^u
 $$
 of this representation space with the group $U(n)^u$.

 Let $F$ be the mapping defined by
 $$
\begin{array}{l}
F: R^\#(U) \times R^\#(U) \cong U(n)^u \times U(n)^u \rightarrow U(n)^u \cong R^\#(U)  \\
( r,s) \mapsto  r \cdot s^{-1} \\
\end{array}
$$

In these terms, $ R^\#(M^3,S_1, \phi)$ is identified as the image of
the pairs  $(a,b) \in R^\#(H_1,\phi) \times   R^\#(H_2) $ for which
$ F( K \cdot J(a), L(b)) = Id \in U(n)^u$. In other words,
for  the composite mapping
$$
G = F( K \cdot J( \star ), L(\star)): ( \  R^\#(H_1,\phi) \times   R^\#(H_2) \ ) \rightarrow R^\#(U) \cong U(n)^u
$$
the inverse of the point $Id \in U(n)^u$ is precisely the set of representations
$R^\#(M^3,S_1, \phi) \cong R(M^3,S_1,[\phi])$.

Since   the spaces  $R^\#(U), R^\#(H_1,\phi),R^\#(H_2)$
are each  of the type $U(n)^x$ for various $x$'s,  they are also  smooth
manifolds.   So once oriented,  making the mappings $K \cdot J, L$ transverse [by a small perturbation,
if needed],  as mappings
of compact, oriented, smooth manifolds, defines an intersection cycle and class
$$
[cyc] := [intersection ] \ \in \ H_W(R^\#(U),Z)
$$
of $( K \cdot J)(R^\#(H_1,S_1, \phi))$ and $L( R^\#(H_2))$ in $R^\#(U)$,
when  $ W = dim \ R^\#(H_1,\phi) + dim \ R^\#(H_2) - dim \ R^\#(U)$ is non-negative.
This intersection is a well defined cycle up to boundaries. It is defined by
moving  the images of the mappings $K \cdot J,L$ to be transverse and  taking the inverse image of $Id$ under $G$.
As will be seen, this oriented intersection cycle regarded as a class in $R^\#(U)$ has
a homological nature.

\vspace{.2in}
It is claimed that
$$
W = (  dim \ U(n)) \ (\chi(S_1) - \chi(M^3)) = ( dim \ U(n))  \cdot T
$$

By definition,
$$
\begin{array}{l}
T = W /(dim \ U(n)) \\
=  ( dim \ R^\#(H_1,\phi) + dim \ R^\#(H_2) - dim  R^\#(U))/( dim \ U(n)) \\
= ((h_1-g_1) + h_2 -u) = (h_1+h_2-u) - g_1
\end{array}
$$
while by $M^3 = H_1 \cup H_2$ with $H_1 \cap H_2 = U$,
$$
\begin{array}{l}
\chi(S_1) - \chi(M^3) \\
= \chi(S_1) - (\chi(H_1) + \chi(H_2) -\chi(U)) \\
= (1-g_1) - ((1-h_1)+(1-h_2) - (1-u)) = (h_1+h_2 -u) - g_1
\end{array}
$$
Hence, $W$ is as claimed.

\vspace{.2in}
\textbf{Orientations}:

To specify the orientations of the above smooth manifolds, it suffices
to orient the tangent spaces at the trivial representation of the
spaces $R^\#(U), R^\#(H_1),R^\#(H_2)$ and the tangent space
at $\phi$ of $R^\#(H_1,S_1, \phi)$. This last is the kernel
of the restriction of the tangent space of $R^\#(H_1)$ at $\phi$
to the tangent space of $R^\#(S_1)$ at $\phi$, which is surjective.
Hence, by definition, it suffices to orient the tangent spaces
at the identity of $R^\#(U), R^\#(H_1), R^\#(H_2), R^\#(S_1)$.
Since each of the fundamental groups is free, these tangent spaces
are precisely, $H^1(U,R)\otimes u(n), H^1(H_1,R) \otimes u(n), H^1(H_2,R) \otimes u(n),
H^1(S_1,R) \otimes u(n)$, where $u(n)$ denotes the Lie algebra of $U(n)$.

Note the long exact Mayer-Vitoris sequence for reduced real cohomology:
$$
\begin{array}{llll}
0 & \rightarrow   H^1(M,R) & \rightarrow H^1(H_1,R) \oplus H^1(H_2,R) & \rightarrow H^1(U,R) \\
& \rightarrow H^2(M,R) & \rightarrow H^2(H_1,R) \oplus H^2(H_2, R) = 0 \oplus 0 & \rightarrow H^2(U,R) = 0
\end{array}
$$
since $H_1,H_2$ are handlebodies, so have the homotopy type of bouquets of circles
and, in particular, have vanishing second homology groups. By $U$ connected, $\tilde{H}^0(U,Q) = 0$.
Hence,
specifying orientations for $H^2(M,R),H^1(M,R)$ determines compatible
choices of orientations for the triple $H^1(U,R)\otimes u(n), H^1(H_1,R) \otimes u(n), H^1(H_2,R) \otimes u(n)$.
In this way, the choices of orientations of $M^3$ and of $H^2(M,R),H^1(M,R), H^1(S_1,R)$
determine the orientation for the cycle
$$
[cyc] \in \ H_W(R^\#(U),Z).
$$
The choice of orientation for $M^3$ determines an orientation on $S_1 \subset \partial M^3$
and so orientations for each of the boundary circles $\partial S_1 = \partial U$.
Also, the induced orientation from $M^3$ on $H_1$ containing $S_1$ give
an orientation on $U \subset H_1$ which together with the given
orientations of the boundary circles of $\partial U$ specify an compatible orientation
of $H^2(U,\partial U,R)= R$.
By this means the class $[cyc]$ is well defined, given
the four orientations of $M^3$ and of   $H^2(M,R), H^1(M,R), H^1(S_1,R)$ which
are together recorded here by the letter $o$.
Since  $(-1)^{dim \ u(n)}   = (-1)^{n}$, changing
the orientation of any one of these four terms   introduces a change of sign by $(-1)^n$
in the value of the cycle $[cyc]$.

\vspace{.2in}
Alternatively and quite generally,  by facts about cap products and Poincar\'{e} duality, employing the mapping
$$
G = F( K \cdot J( \star ), L(\star)): R^\#(H_1,\phi)  \times  \ R^\#(H_2) \rightarrow R^\#(U),
$$
the intersection class $[cyc] \in H_W(R^\#(U),Z)$ may be described as follows:

Let $\{R^\#(U)\}$ denote the specified oriented generator of the top cohomology, \newline
$H^{dim \ R^\#(U)}(R^\#(U),Z) \cong Z$. Let  $[R^\#(H_1,\phi)], [R^\#(H_2)]$ denote the specified oriented generators
in top dimensional homology of $H_{dim \ R^\#(H_1,\phi)}(R^\#(H_1,\phi), Z ) \cong Z$ and
$H_{dim \ R^\#(H_2)}(R^\#(H_2), Z ) \cong Z$, respectively; then there is the equality of
homology classes
$$
\begin{array}{l}
[cyc] = G_\star (  ( [R^\#(H_1,\phi)] \times  [R^\#(H_2)]) ) \ \cap \  G^\star( \{R^\#(U)\})   \ \
  \in   H_W(R^\#(U),Z) \\
\end{array}
$$
Here $\cap$ is the cap product pairing from (homology, cohomology) to homology.
This demonstrates the intrinsically homological (in the representation spaces) nature of these invariants.

Additionally, any two irreducible representations $\phi's$ can be continuously deformed to  each other
through irreducible representations;
so the homology class $[cyc]$ is left unchanged. Hence, invariants based on $[cyc]$  are  independent of
the choice of $\phi$.

\vspace{.2in}
In the case of equality, $\chi(M^3) = \chi(S_1)$, the associated
cycle is a 0-cycle. Adding up  the intersection points according to their signs, after a generic
perturbation, gives a sum of points with signs and multiplicities
which  is the desired integer, assuming that the individual
pieces, $R^\#(U), R^\#(H_1,S_1, \phi), R^\#(H_2)$ which are all of the
form $U(n)^x$ for various $x's$,  have  been oriented compatibly given
$o$. Let $\hat{U}$ denote the closed Riemann surface obtained from $U$
by collapsing individually each boundary circle to a point. Let $u$
be the genus of $\hat{U}$.
In order to get a result independent of the choice of Heegaard decomposition
an additional sign is needed, $(-1)^{(dim \ U(n) \cdot u} = (-1)^{ n \cdot  u} $,
see \S \ref{sectionstable} where the question of stabilization is explicitly
addressed.

The reader will note that the symbol $o$ is dispensed with in the case $T=0$,
corresponding to the integer invariant. The reason is that for
$\lambda_{U(n)}(M^3,S_1)$, (resp., $\lambda_{SU(N)}(M^3,S_1)$), to not
vanish, necessarily $H_2(M^3,R) = 0$ and $H_1(S_1,R) \rightarrow H_1(M^3,R)$
is an isomorphism [see below]. Thus in that case, the only one where signs
are an issue, one can just chose
$o$ so that the orientation of $H_1(S_1,R)$ and $H_1(M^3,R)$ correspond.
With this orientation convention, the associated numerical invariant
depends only on that of $M^3$, whence the notion $\lambda_{U(n)}(M^3,S_1)$,
(resp., $\lambda_{SU(n)}(M^3,S_1)$).

That is, with these  orientation conventions, the definition
of the required integer-valued invariant  in the case $T = \chi(S_1) - \chi(M^3) = 0$ is:
$$
\begin{array}{l}
 \lambda_{U(n)}(M^3,S_1)  \\
\hspace{.5in}  = (-1)^{ n \cdot u} \ [intersection \ number \
 of \ R^\#(H_1,S_1,\phi) \ and \ R^\#(H_2) \ in \ R^\#(U) ]
 \end{array}
 $$
If one takes  absolute values, then  the choice of orientation classes $o$ of
these manifolds will have no effect on the result. Similar remarks apply to
the class $[cyc]$ generally.

\vspace{.2in}
By the above, for $\chi(S_1) - \chi(M^3) =  T=0$, the integer invariant $\lambda_{U(n)}(M^3,S_1)$ has
absolute value:
$$
\begin{array}{l}
|\lambda_{U(n)}(M^3,S_1)|
=    | ( degree \ of \ the \ mapping
\  G : (\  R^\#(H_1,\phi)  \times  \ R^\#(H_2) \ ) \rightarrow R^\#(U) )|
\end{array}
$$
where $| \star |$ denotes the absolute value.

The arguments of this section, by replacing $U(n)$ by $SU(n)$ and $(-1)^{( dim \ u(n)) \cdot u} = (-1)^{n \cdot u}$
by $(-1)^{ (dim \ su(n)) \cdot u } =(-1)^{(n-1) \cdot u}$ defines $\lambda_{SU(n)}(M^3,S_1)$ for the case
$T=0$ and the cycle $[cyc]$ as well.

\section{Proof of theorems \ref{theoreminv1} (resp., \ref{theoreminv1SU}),
 \ref{theoremfirst} and  \ref{theoremapplication}}  \label{sectioncyl}

Homological methods of extending representations from  finite groups were studied in a classical
short paper of  Gerstenhaber and  Rothaus \cite{GerstenhaberRothaus1}.
As in their work,  the degrees of maps between  products of unitary groups will be utilized.

Assume $\chi(M^3) = \chi(S_1)$, i.e., $g= g_1$,  or equivalently,
$ (h_1+h_2)- g_1 = u$.

In order to compute the degree of the mapping
$$
\begin{array}{l}
G : U(n)^{h_1-g_1} \times U(n)^{h_2} \rightarrow U(n)^u \\
R^\#(M^3,S_1,\phi) \times R^\#(H_2) \rightarrow R^\#(U), \\
\end{array}
$$
 it is helpful to review facts about
the cohomology of $U(n)$. A computation of this degree will
prove theorem \ref{theoreminv1} parts 2,3,4.
This detailed discussion of the $U(n)$
case will be followed by discussion of the
corresponding $SU(n)$ case.

 Note that the required mapping
$$
G : U(n)^N \rightarrow U(n)^N,
$$
with $N= (h_1 +h_2) -g_1 = u$ arises in  a very special manner.
There is a homomorphism, say $\rho : F_N \rightarrow F_N$ of free groups
and under the natural identification
$$
Hom(F_N ,U(n)) \cong U(n)^N,
$$
the induced mapping is precisely $G$. Here $Hom(F_N,U(n))$
is the space of
 homomorphisms of $F_N$ to $U(n)$. It is identified with $U(n)^N$
 by sending
$\rho$ to  its evaluation, $ (\rho(y_1), \rho(y_2), \cdots, \rho(y_N))$, for
free generators $y_1,\cdots, y_N$ of $F_N$.

This special character is clear for the mappings $J : R^\#(H_1) \cong U(n)^{h_1} \rightarrow R^\#(U) \cong U(n)^u$
and $K :  R^\#(H_1) \cong U(n)^{h_1} \rightarrow R^\#(U) \cong U(n)^u$ since they are induced by homomorphisms
of the free groups $\pi_1(U) \subset \pi_1(H_1),\pi_1(U) \subset \pi_1(H_2)$. It also holds with the
prescribed identifications for $R^\#(H_1,\phi) \subset R^\#(U)$. Hence, equally so for the composite
$G = F( I \cdot J(\star), K(\star))$, as $F$ is defined by $(r,s) \mapsto r \cdot s^{-1}$.

\vspace{.2in}

The Lie group $U(n)$ is known to have integral cohomology an exterior algebra on generators
$x[j] \in H^{2j+1}(U(n),Z)$ for $ j = 0, \cdots, n-1$, see A. Borel \cite{Borel}. Since $U(n)$ is a group, this cohomology
has a Hopf algebra structure. As is well known, the classes $x[j]$ are primitive, that is, $m^\star (x[j])
= x[j] \otimes 1 + 1 \otimes x[j]$ where $m : U(n) \times U(n) \rightarrow U(n)$ is the group
multiplication mapping, $m(a,b) = a \cdot b$. In particular, the generator of the top
cohomology group $H^{dim \ U(n)}(U(n), Z) = Z$ is represented by the product
$\prod_{j=0}^{n-1} \ x[j]$.

Correspondingly, the product of $N$ copies of $U(n)$ has cohomology
$H^\star(U(n)^N,Z)$ equal to the exterior algebra on generators
$\{ \pi_k^\star (x[j]) \ | \ j=0,\cdots  (n-1), k=1 \cdots N\}$ where
$\pi_k : U(n)^N \rightarrow U(n)$ is the projection onto the $k^{th}$ factor.
Hence, to compute the degree of the mapping $G : U(n)^N \rightarrow U(n)^N$
it will suffice to find the pull backs of the $\pi_k^\star(x[j])$
in $H^{2j+1}(U(n)^N,Z)$; that is, the pull backs of $x[j] \in H^{2j+1}(U(n),Z)$
under the mapping
$$
\pi_k \cdot G : U(n)^N \rightarrow U(n)
$$
where the last sends $(A[1],A[2],\cdots, A[N]) \mapsto w[k]$
where $w[k]$ is a finite product of these matrices to various powers.

For example, the mapping $f : U(n) \rightarrow U(n)$ sending $ A \mapsto A^r$
will, by primitivity, send $x[j]$ to $r \cdot x[j]$ and hence have degree
$r^{n}$ as  a mapping of $U(n) \rightarrow U(n)$.

More generally, for  $f : U(n)^N \rightarrow U(n)$ obtained by
taking products of the matrices \newline $(A[1],\cdots, A[N])$ and their  inverses, if
$A[k]$ appears in the word $w$ in the $A[i]'s$ with total multiplicity $m[k]$,
the sum of powers of its occurrences, then
$$
f^\star( x[j]) = \Sigma_{k=1}^N \ m[k] \ \pi_k^\star (x[j]).
$$

Even more generally, let  $G : U(n)^N \rightarrow U(n)^N$
be the mapping induced by a homomorphism of free groups
$f : F_N \rightarrow F_N$ with $F_N$ a free group
on generators $y_1,\cdots, y_N$,  under the identification
$$
Hom(F_N,U(n)) \cong U(n)^N.
$$
sends the entry $(A[1], \cdots, A[N]) \in U(n)^N$
to $(W[1],W[2],\cdots, W[N])$ with each  $W[j]$ a product
of the matrices $\{A[i]\}$ and their inverses. Let
$m[i,k]$ denote the sum of powers to which $A[k]$
is raised in toto in the product giving the entry $W[i]$.
 By the above
$$
G^\star( \pi_i^\star x[j]) = \Sigma_{k=1}^N \ m[i,k] \ \pi_k^\star( x[j])
$$
independent of $j =0, \cdots (n-1)$.

That is, the primitive classes, the span of the $\pi_k^\star( x[j]), k=1, \cdots, N$,  in dimension $2j+1$,
each transform via the matrix $M = \{ m[i,k]\}$. Since the
generator of $H^{N \cdot( dim \ U(n))}(U(n)^N,Z) = Z$ is precisely
the wedge product $\Lambda_{k=1,\cdots N,j= 0,\cdots ( n-1)} \ \pi_i^\star (x[j])$, the
absolute value of the degree of
the mapping $G$ is just $ |det(M)|^{n} $.
Here the transpose of $M$ is the induced mappings of the
abelianization of $f$,
$
f_\star : H_1(F_N,Z) = Z[y_i, i=1 \cdots N] \rightarrow H_1(F_N,Z) = Z[y_i, i=1 \cdots N].
$ Hence, $
|det(G)| = |( det \ f_\star)|^n .
$

For the example at hand, the abelianization may be described as follows:
The inclusions $S_1 \subset H_1, U \subset H_1, U \subset H_2$ induce maps
of integral cohomologies, each a free abelian group:
$$
\begin{array}{l}
a: H^1(H_1,Z) \rightarrow H^1(S_1,Z), \
 b: H^1(H_1,Z) \rightarrow H^1(U,Z), \ c: H^1(H_2 ,Z) \rightarrow H_1(U,Z)
\end{array}
$$
Let $L$ be the kernel of the mapping $a: H^1(H_1,Z) \rightarrow H^1(S_1,Z)$. Since $H_1$
is obtained from $S_1 \times [0,1]$ by adding additional 1-handles, $\pi_1(H_1) = \pi_1(S_1) \star F_{h_1-g_1}$,
a free product of free groups.
Hence,  $L$ is a direct summand of $H^1(H_1,Z)$ and this mapping is onto. That is,
$H^1(H_1,Z) \cong H^1(S_1,Z) \oplus L$.  Let the composite mapping of free abelian
goups
$$
h : L \oplus H^1(H_2,Z) \subset H^1(H_1,Z) \otimes H^1(H_2,Z) \stackrel{b - c}{\rightarrow} H^1(U,Z)
$$
be called $h$ with  $(b-c)(r,s) = b(r)-c(s)$.
 Then $h$ is the dual of the abelianization of the mapping
of free groups $F_N \rightarrow F_N$   defining the desired mapping $G : U(n)^N \rightarrow U(n)^N$
with $N = (h_1 +h_2) - g_1 = u$ in the example at hand.

In view of the above algebra, the degree of $G$  is zero unless the mapping
$h$ tensor the rationals, a mapping $Q^N \rightarrow Q^N$, is an isomorphism.
In the case that the mapping tensor the rationals is an isomorphism,
the absolute value of the degree of $G$ is
$|H^1(U,Z)/h[L \oplus H^1(H_2,Z)]|^{n}$,
the $n^{th}$ power of the order of the finite abelian group $H^1(U,Z)/h[L \oplus H^1(H_2,Z)]$.

As above, the Mayer Vietoris sequence of the triple $(M^3,H_1,H_2)$ with $U = H_1 \cap H_2$ gives a long exact
sequence in reduced rational cohomology:
$$
\begin{array}{llll}
0 & \rightarrow   H^1(M,Q) & \rightarrow H^1(H_1,Q) \oplus H^1(H_2,Q) & \rightarrow H^1(U,Q) \\
& \rightarrow H^2(M,Q) & \rightarrow H^2(H_1,Q) \oplus H^2(H_2, Q) = 0 \oplus 0 & \rightarrow H^2(U,Q) = 0
\end{array}
$$

Since the image by $h \otimes Q$
of $(L \otimes Q) \oplus H^1(H_2,Q)$ in $H^1(U,Q)$ is contained in  the image of $H^1(H_1,Q) \oplus H^1(H_2,Q)$
in $H^1(U,Q)$, the map induced by $h$ tensor the rationals necessarily is not onto
if $H^2(M^3,Q)$ is non-zero.  Hence, $H^2(M^3,Q) \neq 0$ implies that $\lambda_{U(n)}(M^3,S_1) = 0$.

If $H^2(M^3,Q) = 0$, then $\chi(M^3) = 1 - dim \ H^1(M^3,Q) = 1-g = 1-g_1$, so $H^1(M^3,Q) \rightarrow H^1(S_1,Q)$
is a map of vector spaces of the same rank $g$. If an non-zero  element of $H^1(M^3,Q)$ were in the kernel of this
mapping, then its image lies in $L \otimes Q $ and also maps to a non-zero element of $(L \otimes Q) \oplus H^1(H_2,Q)$
which goes to zero in $H^1(U,Q)$. Hence, unless $H^1(M^3,Q) \rightarrow H^1(S_1,Q)$ is an isomorphism
it follows that $h$ tensor the rationals is not an isomorphism and so the degree of $h$ is zero
and $\lambda_{U(n)}(M^3,S_1) = 0$ again. This gives the proof of theorem \ref{theoreminv1}, part 3.

\vspace{.2in}

In the  case, $H^1(M^3,Q) \rightarrow H^1(S_1,Q)$ is an isomorphism,
one proves $H^2(M^3,Q) =   0$ as follows: By Poincar\'{e} duality
$H^2(M^3,Q) \cong H_1(M^3,\partial M,Q)$ which fits into the long
exact sequence on homology $ H_1(\partial M,Q) \rightarrow H_1(M^3,Q)
\rightarrow H_1(M^3, \partial M,Q) \rightarrow H_0(\partial M,Q) \rightarrow H_0(M^3, Q)$.
The last map is an isomorphism by $M^3, \partial M^3$ connected. Since
$i_\star : H_1(S_1,Q) \rightarrow H_1(M^3,Q)$ is an isomorphism,
a fortiori, $H_1( \partial M,Q) \rightarrow H_1(M^3,Q)$ is onto.
Consequently, by this long exact sequence $H_1(M^3,\partial M,Q) \cong H_2(M^3,Q)^\star$
vanishes.

\vspace{.2in}

There is again
a long exact sequence for reduced integral cohomology:
$$
\begin{array}{llll}
0 & \rightarrow   H^1(M,Z) & \stackrel{i}{\rightarrow}  H^1(H_1,Z) \oplus H^1(H_2,Z) & \stackrel{b-c}{\rightarrow} H^1(U,Z) \\
& \rightarrow H^2(M,Z)  & \rightarrow H^2(H_1,Z) \oplus H^2(H_2, Z) = 0 \oplus 0 & \rightarrow H^2(U,Z) = 0
\end{array}
$$
By exactness and the above, $H^1(M^3,Z) \cong Z^g $ and $H^1(S_1,Z) \cong Z^g$ with the induced
mapping $Z^g \rightarrow Z^g$ rationally an isomorphism. Let $B$ be the kernel of the
 mapping $H^1(U,Z) \rightarrow H^2(M^3,Z)$. The  abelian group $H^2(M^3,Z)$ by the above is of
finite order. Since $h$ factors as $ L \oplus H^1(H_2,Z) \rightarrow B \subset H^1(U,Z)$,
the absolute values of its determinant is the product of the  order of the finite abelian group $H^2(M^3,Z)$ and the
absolute value of the determinant of the truncated mapping
$$
\hat{h} : L \oplus H^1(H_2,Z) \stackrel{j}{\subset} H^1(H_1,Z) \oplus H^1(H_2,Z) \stackrel{b-c}{ \rightarrow} B
$$
of free abelian groups which is an isomorphism when tensored with the rationals.  The required determinant
is just the order of $B/[Image \ \hat{h}]$ times the order $|H^2(M^3,Z)|$.

Additionally, there is the short exact sequence
$0 \rightarrow L  \stackrel{k}{\rightarrow} H^1(H_1,Z)   \stackrel{a}{\rightarrow} H^1(S_1,Z) \rightarrow 0$
which fits together to form the two short exact sequences in the tableau: with $m = a \oplus 0$.
$$
\begin{array}{lllll}
 & & 0 & & \\
 & & \downarrow & & \\
 & & H^1(M,Z) & & \\
 && \downarrow i && \\
 0 \rightarrow    L \oplus H^1(H_2,Z) & \stackrel{k \oplus id }{\rightarrow}  & H^1(H_1,Z) \oplus H^1(H_2,Z) &
 \stackrel{m }{\rightarrow} & H^1(S_1,Z)  \rightarrow 0 \\
 & & \downarrow  (b-c) & & \\
 & & B  & &  \\
 & & \downarrow   &  & \\
 & & 0 & & \\
 \end{array}
 $$
  Define mappings $R: H^1(S_1,Z)/[Image \  m \cdot i] \rightarrow B/[Image \ (b-c) ( k \oplus id)]$ and
 $S :B/[Image \  (b-c) ( k \oplus Id)] \rightarrow   H^1(S_1,Z)/[Image \ m \cdot i] $ by $R(m(x)) = (b-c)(x)$ and
 $S((b-c)(y))) = m(y)$.
 By  exactness, they are well defined and are inverses.

Hence, $|det(h)| = |H^2(M^3,Z)| \cdot | det(\hat{h})| = |H^2(M^3,Z)| \cdot |H^1(S_1,Z)/(m \cdot i)[H^1(M^3,Z)]|$ where $m \cdot i$ appears in the long exact sequence
of the pair $(M,S_1)$:
 $$
 H^1(M^3,Z) \stackrel{m\cdot  i }{\rightarrow} H^1(S_1,Z) \rightarrow H^2(M^3,S_1,Z) \rightarrow H^2(M^3,Z)
\rightarrow H^2(S_1,Z) = 0
$$

 Now the last exact sequence shows that \newline
$|H^2(M^3,S_1,Z)| = |H^2(M^3,Z)| \cdot |H^1(S_1,Z)/(m \cdot  i) [H^1(M^3,Z)]|$ so $det(h) = |H^2(M^3,S_1,Z)|$, as desired.
This then proves that $|\lambda_{U(n)}(M^3,S_1)| = K^{n}$  with
$K = |H^2(M^3,S_1,Z)|$, as claimed.

By definition if there were no  representations extending to $\pi_1(M^3)$
the given representation $\phi$, then  the invariant $\lambda_{U(n)}(M^3,S_1)$
would be zero.
This   observation completes the proof of theorem
\ref{theoreminv1}.

The $SU(n)$ case proceeds in a parallel fashion with minor changes, proving
theorem \ref{theoreminv1SU}.

\vspace{.2in}
\textbf{Proof of theorem \ref{theoremfirst}:}
In the case that $i_\star : H_1( S,Q) \rightarrow H_1(M^3,Q)$
is an isomorphism, by the above argument $H_2(M^3,Q) =0$
so $i_\star : H_j(S,Q) \rightarrow H_j(M^3,Q)$ is an isomorphism for
$j=0,1,2,3$. In particular, $\chi(S) = \chi(M^3)$.

Hence, theorem \ref{theoreminv1} applies to $ (M^3,S)$ to extend
any irreducible representation into $U(n)$.
Since any irreducible $U(n)$ representation
is a sum of irreducibles ones, assembling these
extensions proves theorem \ref{theoremfirst}
in the $U(n)$ case. The $SU(n)$ case
proceeds similarly using theorem \ref{theoreminv1SU}
instead.

\vspace{.2in}
\textbf{Proof of theorem \ref{theoremapplication}:}

To show that $H_2(N_2,Q) \rightarrow H_2(W^3,Q)$ is a rational isomorphism
consider the long exact sequence for the pair $(W,N_1)$. By Poincar\'{e}
duality $0 = H_0(W,N_2,Q) \cong H_3(W,N_1,Q)^\star$, so by the exact sequence
$H_3(W,N_1,Q) \rightarrow H_2( N_1,Q) \rightarrow H_2(W,Q)$, the induced
map $i_\star : H_2(N_1,Q) \rightarrow H_2(W,Q)$ is one to one. Since $H_2(N_1,Q) = Q$
to show that this mapping is an isomorphism it suffices to show that $H_2(W,Q)$
is of rank one. Now since  $H_1(N_1,Q) \rightarrow H_1(W,Q)$ is an isomorphism,
a fortiori, $H_1(\partial W, Q) \rightarrow H_1(W,Q)$ is onto. But by the
long exact sequence of the pair $(W,\partial W)$,
$H_1(\partial W, Q) \rightarrow H_1(W,Q) \rightarrow H_1(W,N_1,Q) \rightarrow
H_0(\partial W,Q) \rightarrow H_0(W,Q)$, this surjection implies that
$H_1(W,N_1,Q)$ maps isomorphically to the kernel of
$H_0(\partial N,Q) \rightarrow H_0(W,Q)$ which has rank one. Again, by Poincar\"{e} duality
$H_2(W,Q) \cong H_1(W,N_1,Q)^\star \cong Q $. In toto,
this proves $i_\star : H_2(N_1,Q) \cong H_2(W,Q)$. That is, $W$ is
a rational homology cobordism.

Take a properly embedded path $\gamma$ from a point, say $p$, of the boundary component  $N_1$ to a point $q$ of
the other boundary component $N_2$ of  $W^3$.
Let $M^3$ be the result of deleting a small tubular neighborhood of the path  $\gamma$. It has one boundary
component, the connected sum of $N_1$ and $N_2$. Suppose that
its boundary, $\partial M^3$, is of genus $g$.
 Let this tubular neighborhood intersect $N_1$ in a 2-disk $D^2$ about $p$.

 If $ S$ is not all of $N_1$,
chose $\gamma$ so that $D^2$ is disjoint from $S \subset N_1$. The manifold $W^3$ can be recovered from
$M^3$ by adding in the tubular neighborhood of $\gamma$. Note that $W^3$ has the homotopy type
of $M^3$ with the 2-disk $D^2$ added, $ W^3 \cong M^3 \cup D^2$.

Note that by assumption,
the restriction mapping $i :H^2(W,Q) \rightarrow H^2(N_1,Q) = Q$ is an isomorphism.

Let $N_1' $ be the closure of $N_1 - D^2$.
Now  $W^3$ has $M^3 \cup D^2$ as a deformation retract which keeps $M$ fixed. Now  identify $N_1$ with $N_1' \cup D^2$, then
 the Mayer-Vietoris sequence gives
the commutative diagram with exact rows for rational coefficients:
$$
\begin{array}{cccc}
H^1(\partial D^2) = Q & \rightarrow H^2(W)) =  H^2(M^3 \cup D^2) & \rightarrow H^2(M) \oplus H^2(D^2)
= H^2(M)   & \rightarrow H^2(\partial D^2) = 0 \\
\downarrow \cong  & \downarrow \cong i   & \downarrow & \downarrow \cong  \\
H^1(\partial D^2) = Q & \stackrel{\cong}{\rightarrow} H^2(N_1) = H^2(N_1' \cup D^2) &
 \rightarrow H^2(N_1') \oplus H^2(D^2) = 0   & \rightarrow H^2(\partial D^2) = 0
\end{array}
$$
This shows that  $ H^2(M^3,Q) = 0$.

Similarly, the restriction mapping $H^1(M^3,Q) \rightarrow H^1(N_1',Q)$ ia an isomorphism from the commutative
diagram of exact sequences for reduced rational cohomology:
$$
\begin{array}{cccc}
0 & \rightarrow H^1(W)) =  H^1(M^3 \cup D^2) & \rightarrow H^1(M) \oplus H^1(D^2) = (H^1(M) \oplus 0)
 & \rightarrow H^1(\partial D^2) = Q \\
\downarrow  & \downarrow \cong i   & \downarrow & \downarrow \cong  \\
0& \rightarrow H^1(N_1) = H^1(N_1' \cup D^2) & \stackrel{\cong}{\rightarrow}  H^1(N_1') \oplus H^1(D^2)
& \stackrel{0}{\rightarrow} H^1(\partial D^2) = Q
\end{array}
$$

Now let $\phi : \pi_2(N_1) \rightarrow U(n)$ be any representation. Since any such is a direct sum
of irreducibles, we may assume $\phi$ is irreducible. Then theorem \ref{theoreminv1} applies
to the manifold $M^3$ and subspace $N_1' \subset \partial M^3$, for the restriction, say $\phi'$ of $\phi$ to $\pi_1(N_1')$.
By that theorem there is a representation $\Phi'$ of $\pi_1(M^3)$ to $U(n)$
which extends $\phi'$ to all of $\pi_1(M^3)$. But since $\phi'(\partial D^2) = \phi(\partial D^2) = Id$,
the representation $\Phi'$ of $\pi_1(M^3)$ extends to a representation, say $\Phi$ of
$\pi_1(W^3) = \pi_1(M^3 \cup D^2)$ which restricts precisely to $\phi$.

\vspace{.2in}
The $SU(n)$ case proceeds in a parallel fashion.
The integral cohomology of $SU(n)$ is a Hopf-algebra
generated by primitive classes in dimensions $3,5,7,\cdots, 2n-1$
with cup product the generator of the top cohomology $H^{n^2-1}(SU(n),Z)$.
For this reason, for $T=0$ the associated degree in this $SU(n)$
case, computing the invariant $\lambda_{SU(n)}(M^3,S_1)$, is
$K^{n-1}$ in view of the $n-1$ primitive generators. In a similar
fashion, the dependence on the change of orientation goes as
$  (-1)^{n-1}$ and the definition of the invariant
has the stabilizing term $(-1)^{(n-1) \cdot u}$ with
$u$ the genus of $\hat{U}$, the closed Riemann surface obtained
from $U$ by collapsing each individual boundary circle of $U$ to a point.

\section{Stabilization and invariants for $T >0$.}  \label{sectionstable}

The classical results of J. Singer \cite{Singer1} imply that any two handlebody
decompositions adapted as here to the boundary decomposition
$\partial M_3 = S_1 \cup S_2$ are related by isotopies and
a standard stabilization process. This applies equally well to
the special decomposition considered here which arise by
adding one handles to the regular neighborhoods of $S_1, S_2 $ in $\partial M^3$.
The stabilization process has an explicit description: Given $H_1,H_2$ one chooses
a 3-ball, say $D^3$ in the interior of $M^3$ which intersects
$U$ in a 2-disk, say $D^2$ and $H_1,H_2$ in two 3-balls, say
$B_1,B_2$ respectively. Then this decomposition into $H_1,H_2$
is modified by replacing $D^3 \cap H_1$ by $B_1$ union an
additional handle obtained by connecting two interior
points of $D^2$ by and an arc in $B_2$ and thickening it.
If this new handlebody is called $H_1'$ and the
closure of $M^3 - H_1'$ is called $ H_2'$, then
$(H_1',H_2')$  form a handlebody decomposition
adapted to the pair $(S_1,S_2)$. Let $U' = H_1' \cap H_2'$.
The genus of $U'$ is that of $U$ increased by 1.

Since the union $D^3 \cup H_j'$ is isotopic to $H_j$, j=1,2,   and
the union $U' \cup D^3$ has $U$ as a deformation retract, the
inclusions $H_j' \cup D^3 \cup H'_j \cong H_j$ and $U' \subset U' \cup D^3 \sim U$ define standard
inclusions:
$$
R^\#(H_1) \subset R^\#(H_1), \ R^\#(H_2) \subset R^\#(H_2), R^\#(U) \subset R^\#(U'), R^\#(H_1,\phi) \subset R^\#(H'_1,\phi)
$$
A simple check shows that the pointwise intersection
$$
  [ Image \ R^\#(H_1,S_1,\phi) ] \cap [ Image \ R^\#(H_2) ] \ in \ R^\#(U)
 $$
 equals the intersection
 $$
 [ Image \ R^\#(H'_1,S_1,\phi) ] \cap [ Image \ R^\#(H'_2) ] \ in \ R^\#(U')
$$

Correspondingly, making these transverse in $U$ and extending this to $U'$ shows
that the intersection cycle $[cyc] \in H_W(R^\#(U),Z)$  maps under the inclusion
$R^\#(U) \subset R^\#(U')$ to the intersection cycle $[cyc'] \in H_W(R^\#(U'),Z)$
However, a close check of signs shows that with the standard orientations
for the choice $o = orientation \ of  \ (H^1(M,R), H^2(M,R), H_1(S_1,R))$ the stabilization
process gives a sign change of $(-1)^{dim \ U(n)}$ under the process of one stabilization.

Again, let $\hat{U}$ be the closed surface obtained by collapsing each component of the
boundary of $U$ individually to a point and $u$ denote the genus of $\hat{U}$.

\vspace{.2in}
Hence, defining
$$
proper[cyc] = (-1)^{ ( dim \ U(n)) \cdot u} \ [cyc] \in H_W(R^\#(U),Z)
$$
yields a cycle which is natural for the stabilization process.
This cycle is in dimension
$$
W = ( dim \ U(n)) T = ( dim \ U(n)) \ ( \chi(S_1) - \chi(M^3)).
$$

\vspace{.2in}
In order to have universal cohomology classes  on which to evaluate
these cycles, one turns to the work of Atiyah and Bott
on Yang-Mills gauge theory in the context of Riemann surfaces \cite{AtiyahBott1}.
There they consider the space of $U(n)$ representations,
$$
R_{U(n)}^\#(\Sigma) = Hom( \pi_1(\Sigma), U(n)),
$$
of the fundamental group of a closed Riemann surface of genus $g$
and compute the rational homology
of the quotient by conjugation,
$$
\mathcal{R}_{U(n)}(\Sigma) = R_{U(n)}^\#(\Sigma)/ conjugation.
$$

The classifying space for flat $U(n)$ connections
is the classifying space $B \mathcal{G}_{U(n)}$,
with $\mathcal{G}_{U(n)}$ the gauge group of
continuous mappings,
$$
\mathcal{G}_{U(n)} = Map(\Sigma, U(n));
$$
so there are natural mappings
$$
\mathcal{R}_{U(n)}(\Sigma) \rightarrow B Map(\Sigma, U(n)) = B \mathcal{G}_{U(n)}
$$
which stabilize well under increasing genus.

As it happens, there are universal classes for $H^\star(B Map(\Sigma,U(n)),Q)$ :
$2g$  exterior generators each in dimensions $1,3,\cdots, 2n-1$
[corresponding to the first homology];
one polynomial generator each in dimensions $2,4,\cdots, 2n$
and one polynomial generator in dimensions $2,4,\cdots, 2n-2$ [ corresponding
to second homology]; one polynomial generator each in dimensions
$2,4,\cdots, 2n$ [corresponding to the zero-th homology].
Altogether this is a product of free exterior and polynomial
generators. The generators are labeled by the chosen basis
for $H^1(\Sigma,Z), H_2(\Sigma,Z) = Z, H^0(\Sigma,Z) = Z$.

Under increasing genus these generators map by restriction
to corresponding generators [or to zero if that part of the homology
dies in the restriction].  Hence, one may potentially     form universal
classes which when evaluated
on the proper cycles,
$
[cyc]
$,
give invariants independent of the handlebody
chosen in computing them.

Naturally, there are corresponding
statements for the case $SU(n)$.

\vspace{.2in}
Since Atiyah and Bott's work \cite{AtiyahBott1} dealt with representations of
fundamental groups of closed Riemann surfaces, we must
restrict the chosen representation $\phi : \pi_1(S_1,p) \rightarrow U(n)$
to have the property that it sends each boundary component of the
subsurface $S_1$ to the identity. Then the intersection cycle will
lie in the smaller space of representations which are restrictions
of the fundamental group of $\hat{S}_1$, the closed Riemann surface
obtained by filling in the boundary disks.

In this manner, for each of Atiyah and Bott's universal
characteristic classes \cite{AtiyahBott1}, there is an associated
invariant. This is parallel to the 4-manifold  invariants
of Donaldson \cite{Donaldson1}.

\vspace{.2in}

Two significant issues remain in the $U(n)$ case. The cycles defined
here lie in $R^\#(\hat{U})$
 for the associated closed surface $\hat{U}$
 once the above restriction on $\phi$ is adopted.
However, in the paper of Atiyah and Bott \cite{AtiyahBott1} the
 classes factor through $\mathcal{R}(\hat{U})$.
 Hence, in the universal example, the
 gauge group $\mathcal{G}_{U(n)}$  is to be replaced by the smaller
  group $\mathcal{G}_{U(n)}'  = Map'( \Sigma, U(n))$ of maps
  which are pointed, i.e., mappings which send a fixed
  base point to the identity. Atiyah and Bott show that
  the integral cohomology of $B \mathcal{G}_{U(n)}$
  is the tensor product of the integral cohomology of $B U(n)$
  with the integral cohomology of $\mathcal{G}_{U(n)}' $.

  The cohomology of the classifying space
  $B \mathcal{G}_{U(n)}'  $ lacks the
  generators one each in dimensions
  $2,4,\cdots,2n$ corresponding to
  the zero-dimensional generator
  $H^0(\Sigma,Z) = Z$.

\vspace{.2in}

Secondly, since the mapping to the universal space utilizes
a choice of basis for the homology $H^1(\hat{U},Z)$
to identify the generators, in order to produce
 true invariants independent of these choices,
it is necessary to look for elements of
the algebra invariant under the automorphisms
of $H^1(\Sigma,Z)$ endowed with its cup pairing,
i.e., invariant under the natural action of
the discrete symplectic group $Sp(2g,Z)$.

Consequently, each of the 2g exterior generators
in dimension $2k-1$  contribute but one invariant class,
while all the polynomial generators survive. Passing to the
limit, the universal symplectically invariant classes
form an algebra with:
polynomial generators in  dimensions
$2,6,10, 2(2n-1)$ [from the
symplectic form on the odd generators], two polynomial
generators in each even dimension [ from the polynomial
classes corresponding to $H^2(\Sigma,Z)$.

These classes stabilize to give
two polynomial generators
$x_i,y_i$ with $x_i$ in dimension
$2i$ and $y_j$ in dimension $4j-2$, the
second from the symplectically invariant
classes.

\vspace{.2in}
In consequence, for any monomial
$$
X _{I,J}= ( x_{i_1}^{r_1} x_{i_2}^{r_2} \cdots x_{i_a}^{r_a})  ( y_{j_1}^{s_1} y_{j_2}^{s_2} \cdots y_{j_b}^{s_b})
 $$
 with
 $$
 T = \chi(S_1)- \chi(M^3) = (\  \Sigma_{p=1}^a (2  i_p)(r_p)  \ + \ \Sigma_{p=1}^b ( 4  j_p -2)(s_p)  \ )
$$
 there is an invariant $\lambda_{U(n),I,J}( M^3, S_1, \phi)$  counting with signs and multiplicities the
number of  $(dim \ U(n)) \times T$ parameter families of representations which extend to $M^3$
the irreducible representation $\phi$, where $\phi$ is chosen to send each boundary component
of $S_1$ to $Id$.

These may be codified as a single homogeneous polynomial
$$
\Lambda_{U(n)}(M,S_1,o) = \Sigma_{I,J} \ \lambda_{U(n),I,J}(M,S_1,o) \ X_{I,J}
$$
in a manner reminiscent of the work of Donaldson on 4-manifolds \cite{Donaldson1}.

\vspace{.2in}
For the $SU(n)$ case, the cohomology of the classifying
space for
  $B \mathcal{G}_{SU(n)}$ is a exterior algebra on
  classes in dimensions $3,5,\cdots, 2n-1$ and
  polynomial algebras on generators in dimensions
  $4,6,\cdots, 2n$ and also $ 4,6,\cdots, 2n-2$.

   Hence, deleting the polynomial invariant arising from $H^0(\hat{U},Z)$
   and replacing the exterior generators by the symplectically invariant
   ones, we arrive at polynomial generators $x_i$ in dimension $2i$ for
   $ i>1$ and polynomial generators $y_i$ in dimensions $4j-2$ for $ j >1$.

   The corresponding multi-index invariants for $SU(n)$, $\lambda_{SU(n),I,J}(M^3,S_1,o)$,  are indexed by multi-indices
(possibly vacuous) $ I = [(i_1,r_1), (i_2,r_2), \cdots, (i_a,r_a)], \ J =  [(j_1,s_1), (j_2,s_2), \cdots, (j_b,s_b)] $ of non-negative
indices with  $i_1 < i_2< \cdots i_a$ and $ j_1 < j_2 < \cdots < j_b$ subject to the added constraints
$  i_1>1$ if $I$ is present  and $j_1 >1 $ if  $J$  is present. Evaluating these on $[cyc]$
yields the desired $SU(n)$ invariants:
For
$$
 T = \chi(S_1)- \chi(M^3) = (\  \Sigma_{p=1}^a (2  i_p)(r_p)  \ + \ \Sigma_{p=1}^b ( 4  j_p -2)(s_p)  \ )
$$
 there is an invariant $\lambda_{SU(n),I,J}( M^3, S_1, \phi)$  counting with signs and multiplicities the
number of  $(dim \ SU(n)) \times T$ parameter families of representations which extend to $M^3$
the irreducible representation $\phi$, where $\phi$ is chosen to send each boundary component
of $S_1$ to $Id$.

All these may be recorded by a single homogeneous polynomial invariant
$$
\Lambda_{SU(n)}(M,S_1,o) = \Sigma_{I,J} \ \lambda_{SU(n),I,J}(M,S_1,o) \ X_{I,J}
$$
with the sum over multi-indices $I,J$ as above with the added constraints $i_1 >1$ if $I$ present and if $j_1>1$
if $J$ is present.

  With these choices theorems \ref{theoreminv2}, \ref{theoreminv2SU} holds.

\vspace{.2in}

\begin{example} Non-trivial homogeneous polynomial invariants. \end{example}
A simple example in which these polynomial invariants are non-vanishing is provided as follows:

Let $S$ be a Riemann surface of genus $g$  with one boundary component and $M^3 = S \times [0,1]$.
Now choose a simple closed curve, say $\gamma$, which separates $S \times 0$ into a
genus $h$ subsurface, say $S_1$, with single boundary component $\gamma$ and remaining genus $g-h$ subsurface,
say $T$. Now suppose that the genus of $S_1$ is at least two, and an irreducible representation
$\phi : \pi_1(S_1) \rightarrow U(n)$ is chosen, as may be done, so that $\phi$ sends the
boundary component $\gamma$ to the identity. In the standard way, $\pi_1(S_1)$ is a free group
on generators, say $\{a_j,b_j \ | \ j=1,\cdots, h\}$,  with $\gamma $ represented by the product of
commutators $\prod_{j=1}^h \ [a_j,b_j]$. Then $\pi_1(S \times 0)$ may be presented as
a free group of these $2h$ generators plus $2(g-h)$ more, say $\{a_j,b_j \ | \  j= h+1,\cdots g\}$,
where the boundary of $S \times 0$ is represented by the product of commutators $\prod_{j=1}^g \ [a_j,b_j]$.
Since the fundamental groups of $S \times 0$ and $M^3 = S \times [0,1]$ are isomorphic,
the space of representations of $\pi_1(M^3) $ into $U(n)$ is just the product
$R^\#(M) =  \times_{k=1}^{2g} U(n)$ and the space of of representations into $U(n)$
restricting to $\phi$ is just the product $R^\#(M,S_1,\phi ) = \times_{k=2h+1}^{2g} U(n) = U(n)^{2(g-h)}$.
The cohomology  class  $y_j$ pulls back to this product as a sum of terms $\Sigma_{j=1}^{g-h} \pi_j^\star(a_j \otimes b_j)$
where $\pi_j$ is the projection to the $j^{th}$ pair $U(n) \times U(n)$ in $U(n)^{2(g-h)} = (U(n) \times U(n))^{g-h}$ and $a_j \otimes b_j$
is the class in $H^{2j-1}(U(n),z) \otimes H^{2j-1}(U(n),Z)$ representing the  product of the
primitive classes, here $a_j, b_j$, in these dimensions and factors. The product of the $a_j$'s over $j = 1 ,\cdots, n$ evaluates
the fundamental class of $U(n)$ to one as does the product of the $b_j$'s.
Consequently,  the characteristic class
$( \prod_{j=1}^n \ y_j)^{g-h}$ in dimension $2 (g-h) \cdot ( dim \ U(n))$ evaluates
the cycle  $U(n)^{ 2 ( g-h)}$ giving   $((g-h)!)^n$ up to sign.

Similarly for $SU(n)$,  evaluating
the characteristic class $( \prod_{j=2}^n \ y_j)^{g-h}$ in dimension $ 2(g-h) \cdot ( dim \ SU(n))$
on the corresponding cycle $SU(n)^{h-g}$  yields up to sign $((g-h)!)^{n-1}$.

\vspace{.2in}
Questions about extending various kinds of representations arise naturally  in low dimensional topology and in group theory,
e.g., \cite{GerstenhaberRothaus1,CochranOrrTeichner1,CochranOrrTeichner2,CochranTeichner,CochranKim,Baumslag}.
It is tempting  to wonder if the present results can be profitably generalized to the combinatorial group
theory setting.

\section{Gauge Theoretic Reformulation}  \label{sectiongauge}

It's natural to reformulate the present $U(n)$ (respectively, $SU(n)$) representation-theoretic numerical invariants in terms
of gauge theory of flat $U(n)$ (resp., $SU(n)$) connections, in analogy to the work of Taubes on $SU(2)$ flat
connections \cite{Taubes1}. This can be carried out by utilizing the $Z/2Z$ valued spectral flow
invariant in the setting of real bounded Fredholm operators of index 0 acting on a separable Hilbert space
, see  \cite{Taubes1} page 571 or generally the `real K theory', denoted $KR(X)$, of Atiyah, Patodi, and Singer,
or the direct analysis of Koschorke
\cite{AtiyahPatodiSinger1,Atiyah1,Koschorke1}, in place of the $Z$ valued spectral flow invariant of
Atiyah, Patodi, and Singer \cite{AtiyahPatodiSinger2, BoossLeschPhillips1}
in the setting of self-adjoint operators utilized  by Taubes. This will be discussed in a future paper.

Such a reformulation offers the possibility of defining Floer type homology groups \cite{Floer1, Floer2, Floer3}
in these $U(n)$ and $SU(n)$ settings for manifolds with boundary.
However, carrying
out a gauge-theoretic reformulation for the more general homogeneous polynomial invariants introduced
 here would be more challenging.

\frenchspacing

\bibliographystyle{plain}

\begin{thebibliography}{99}




\bibitem{Akbulut1} S. Akbulut \& J. McCarthy, {Casson's invariant for Oriented Homology Three-Spheres, An Exposition},
Princeton University Press, Mathematical Notes \# 36, 1986.

\bibitem{Atiyah1} M. F. Atiyah, {K Theory and Reality}, Quar. J. Math., 17  (1966),  367-386.


\bibitem{AtiyahBott1} M. F. Atiyah, R. Bott, {Yang Mills Equations over Riemann Surfaces}, Phil. Trans.
R. Soc. Lond.,A \ 308 (1982), 523-615.


\bibitem{AtiyahPatodiSinger1} M. F. Atiyah, V. K. Patodi \& I. M. Singer, {Index of Elliptic Operators V},
   Ann. of Math., Second Series,  93 (1971),  139-149.



    \bibitem{AtiyahPatodiSinger2} M. F. Atiyah, V. K. Patodi \& I. M. Singer,
    {Spectral asymmetry and Riemannian geometry:  III}, Math. Proc. Cambridge Philos. Soc., 79 (1976), 71-99.

\bibitem{Baumslag} G. Baumslag, R. Mikhailov, K.Orr, {A new look at finitely generated
meta-abelian groups}, in
Computational and combinatorial group theory and cryptography,  Contemp. Math., (2012),
 Vol. 582.

\bibitem{BodenHerald2} H. Boden, C. Herald, {A Connected Sum Formula for the SU(3) Casson Invariant },
J. Diff. Geom., 53 (1998), 443-464.


\bibitem{BodenHerald1} H. Boden, C. Herald, {The SU(3) Casson invariant for Integral Homology 3-Spheres},
J. Diff. Geom., 50 (2001), 147-206.



\bibitem{BodenHeraldKirkKlassen1} H. Boden, C. Herald, P. Kirk, E. Klassen, {Gauge theoretical Invariants of Dehn Surgeries
on Knots}, Geom. Topol., 5 (2001), 143-226.



\bibitem{BodenHeraldKirk1} H. Boden, C. Herald, P. Kirk, {The Integer Valued $SU(3)$ Casson
Invariant for Brieskorn Spheres}, J. Diff. Geom., 71 (2005), 23-83.



\bibitem{BodenNicas1} H. Boden and A. Nicas, {Universal Formulas for $SU(n)$ Casson Invariants of Knots},
Trans. Amer. Math. Soc, 352 (2000), 3149-3187.


\bibitem{BoossLeschPhillips1} B. Booss-Bavnbek, M. Lesch, and J . Phillips,
{ Unbounded Fredholm Operators and Spectral Flow}, 57 (2005), 225-250.


\bibitem{BoyerLines1} S. Boyer, D. Lines, { Surgery Formulae for Casson's invariant and Extension to Homotopy Lens Spaces},
J. Rien Angew. , 405 (1990), 181-220.

\bibitem{BoyerNicas1} H. Boyer and A. Nicas {Varieties of Group Representations and Casson's
Invariants for Rational Homology 3-Spheres},
Trans. Amer. Math. Soc, 322 (1990), 507-522.


\bibitem{Borel} A. Borel,  {Sur L'Homologie et la Cohomologie des Groupes de Lie Compacts Connexes},
 American Journal of Mathematics,  76,  \# 2 (1954), 273-342.


\bibitem{CappellLeeMiller1} S. Cappell, R. Lee, E. Miller, {A Perturbative SU(3) Casson Invariant},
Commentarii Mathematici Helvetici, 77 (2002),  491-523

\bibitem{CochranOrrTeichner1} T. Cochran, K. Orr, P. Teichner,   {Knot concordance, Whitney towers and $L^2$ signatures},
Journal-ref: Ann. of Math. (2), Vol. 157 (2003), no. 2, 433-519.


\bibitem{CochranOrrTeichner2} T. Cochran, K. Orr, P. Teichner, {Structure in the classical knot concordance group},
Journal-ref: Comment. Math. Helv. 79 (2004) 105-123.

\bibitem{CochranTeichner}T. Cochran, P. Teichner, {Knot concordance and von Neumann $p$-invariants},
Journal-ref: Duke Math. Journal, 137 (2007), no.2, 337-379.

\bibitem{CochranKim} T. Cochran, T. Kim, {Higher-order Alexander invariants and filtrations of the knot concordance group}.
Journal-ref: Trans. Amer. Math Soc, 360 no.3 (2008), 1407-1441.



\bibitem{Donaldson1} S. Donaldson, { Polynomial Invariants for Smooth Four-Manifolds},
Topology 29 (1990), 257-315.

\bibitem{Floer1} A. Floer, {An instanton Invariant for 3-Manifolds}, Comm. Math. Phys., 118 (1988), 215-240.



\bibitem{Floer2} A. Floer, { Morse Theory for Lagrangian Intersections}, J. Diff. Geom., 28 (1988), 513-547.

\bibitem{Floer3} A. Floer, { Instanton Homology, Surgery, and Knots}, in: Geoemtry of Low-Dimensional Manifolds, 1 (Durham 1989),
97-114, London Math. Soc. Lecture Note Ser., 150, Cambridge Univ. Press, 1990.

\bibitem{Frohman1} C. Frohman,  {Unitary Representations of Knot Groups},
Topology 32 (1993), 121- 144.

\bibitem{FrohmanLong1} C. Frohman and D.  Long, { Casson's Invariant and Surgery on Knots},
Proc. of the Edinburgh Math. Soc., 35 (1992), 383-395.

\bibitem{Goldman1} W. Goldman, {The Symplectic Nature of the Fundmantal Groups of Surfaces}, Adv. in Math.,
54 (1984), 200-225.


\bibitem{Goldman2} W. Goldman, {Representations of fundamental groups of surfaces},
in Geometry and Topology, Proceedings, University of Maryland 1983-
1984, J. Alexander and J. Harer (eds.), Lecture Notes in Math-
ematics, vol. 1167 (1985), 95-117, Springer-Verlag New York.


\bibitem{GerstenhaberRothaus1} M. Gerstenhaber and O. Rothaus, {On the Set of Solutions
of Equations in Groups},  Proc. Nat. Acad. Sci. U.S.A., 49 (1962), 1531-1533.

\bibitem{Herald1} C. Herald,  {Flat Connections, the Alexander Invariant, and Casson's Invariant},
Comm. Anal. Geom. 5 (1997), 93-120.

\bibitem{Koschorke1} U. Koschorke, {Infinite Dimensional K-theory and Characteristic Classes of Fredholm Bundle
Maps}, Proc.Sympos.Pure Math.,Amer.Math.Soc, 17 (1970).

\bibitem{Lescop1} C. Lescop, {Global Surgery Formula for the Casson-Walker Invariant}, Princeton Univ. Press,
1996.

\bibitem{Lin1} X. Lin, {A Knot Invariant via Representation Spaces}, J. Diff. Geom, 35 (1992), 337-357.

\bibitem{Nicas1} A. Nicas,  { Combinatorial Identities in the Theory of $SU(n)$ Casson Invariants
of Knots}, Pure \& Applied Math. Quarterly, 2 ( 2006), 795-816.

\bibitem{RubermanSaveliev1} D. Ruberman and N. Saveliev,  { Casson-Type Invariants in Dimension Four},
Geometry and Topology of Manifolds, Fields Inst. Commun., 47 (2005), 281-306.






\bibitem{Saveliev1} N. Saveliev, {Lectures on the Topology of 3-Manifolds : An Introduction
to the Casson Invariant}, Walter De Guryter, Berlin, New York, 1999.

\bibitem{Saveliev2} N. Saveliev, {Invariants of  Homology 3-Spheres}, Springer-Verlag Berlin, Heidelerg, New York (2002)

\bibitem{Singer1} J. Singer, {Three Dimensional Manifolds and Their Heegaard Diagrams}, Trans. Amer. Math. Soc., 35 (1933), 88-111.

\bibitem{Taubes1}  C. H. Taubes,
{Casson's Invariant and Gauge Theory},
J. Differential Geoometry,
31 (1990),  547-599.

\bibitem{Walker1} K. Walker, {An Extension of Casson's Invariant}, Annals of Mathematical Studies, \# 126, Princeton Press,
Princeton, New Jersey, (1992).


\end{thebibliography}

\vspace{.2in}

Sylvain E. Cappell, Courant Institute of Mathematical Sciences, New York University. cappell@cims.nyu.edu

\vspace{.2in}

Edward Y. Miller, Mathematics Department, Polytechnic School of Engineering,
New York University. emiller@poly.edu

}

\end{document}